\documentclass[10pt]{amsart}
\usepackage{amssymb}

\usepackage{amscd}
\usepackage{amsthm}
\usepackage[dvips]{graphicx}

\usepackage[all]{xy}

\newtheorem{Lemma1}{{Lemma}}
\newtheorem{Theo1}[Lemma1]{{Theorem}}
\newtheorem*{Theo2}{{Theorem}}
\newtheorem{Def1}[Lemma1]{{Definition}}
\newtheorem{Prop1}[Lemma1]{{Proposition}}
\newtheorem{Claim1}[Lemma1]{{Claim}}
\newtheorem{Rem1}[Lemma1]{{Remark}}
\newtheorem{Cor1}[Lemma1]{{Corollary}}
\newtheorem{Ex1}[Lemma1]{{Example}}

\newenvironment{Lemma}{\begin{Lemma1}}{\end{Lemma1}}

\newenvironment{Prop}{\begin{Prop1}}{\end{Prop1}}
\newenvironment{Rem}{\begin{Rem1}\rm}{\end{Rem1}}
\newenvironment{Theorem}{\begin{Theo1}}{\end{Theo1}}

\newenvironment{Cor}{\begin{Cor1}}{\end{Cor1}}

\title{K\"ulshammer ideals of algebras of quaternion type}
\author{Alexander Zimmermann}

\address{Universit\'e de Picardie,
D\'epartement de Math\'ematiques et LAMFA (UMR 7352 du CNRS),
33 rue St Leu, F-80039 Amiens Cedex 1, France}
\email{alexander.zimmermann@u-picardie.fr}
\date{May 24, 2016; revised January and May 2017}

\newcommand{\lra}{\longrightarrow}

\newcommand{\ra}{\rightarrow}
\newcommand{\sdp}{\times\kern-.2em\vrule height1.1ex depth-.05ex}
\newcommand{\epi}{\lra \kern-.8em\ra}

\newcommand{\N}{{\mathbb N}}

\newcommand{\soc}{\text{soc}}
\newcommand{\dickebox}{{\vrule height5pt width5pt depth0pt}}

 \normalbaselines
\subjclass[2000]{Primary: 16G10, 20C05;
Secondary: 18E30, 16G60}
\keywords{derived equivalences, stable equivalences of Morita type, algebras of quaternion type, tame blocks, K\"ulshammer ideals, socle deformation}

\begin{document}

\begin{abstract}
For a symmetric algebra $A$ over a field $K$ of characteristic $p>0$
K\"ulshammer constructed a descending sequence of ideals of the
centre of $A$. If $K$ is perfect this sequence was shown to be an
invariant under derived equivalence and for algebraically closed
$K$ the dimensions of their image in the stable centre was shown to be invariant
under stable equivalence of Morita type.
Erdmann classified algebras of tame representation type which
may be blocks of group algebras, and Holm classified Erdmann's list
up to derived equivalence. In both classifications certain parameters
occur in the classification, and
it was unclear if different parameters lead to different algebras.
Erdmann's algebras fall into three classes,
namely of dihedral, semidihedral and of quaternion type.
In previous joint work with Holm we used K\"ulshammer ideals to distinguish
classes with respect to these parameters in case of algebras of dihedral and
semidihedral type. In the present paper we determine the K\"ulshammer
ideals for algebras of quaternion type and distinguish again algebras with
respect to certain parameters.
\end{abstract}

\maketitle

\section*{Introduction}

Erdmann gave in \cite{Erdmann} a list of basic symmetric
algebras of tame representation type which include all the algebras which
may be Morita equivalent
to blocks of finite groups of tame representation type. She obtained these algebras
by means of properties of the Auslander-Reiten quiver which are
known to hold for blocks of group rings with tame representation type.
These algebras
are subdivided into three classes, those of dihedral type, of
semidihedral type and of quaternion type, corresponding to the
possible defect group in case of group algebras, and actually defined
by the behaviour of their Auslander-Reiten quiver.
Holm refined in \cite{Holmhabil}
Erdmann's classification of those algebras which may occur as blocks of group algebras to a classification up to derived equivalence.
However, in \cite{Erdmann,Holmhabil} the algebras
are defined by quivers with relations, and the relations involve
certain parameters, corresponding mostly to deformations of the socle of the algebras.
It was unclear in some cases if different parameters lead to
different derived equivalence classes of algebras. The question of non trivial socle deformations
appears to be a very subtle one in this special case, but also in general,
and little progress was made on this question until
very recently.

In \cite{kuelsquest} we showed that a certain sequence of ideals of
the centre of a symmetric algebra defined previously by K\"ulshammer \cite{Ku1}
is actually invariant under derived equivalences if the base field is perfect. We call this
sequence of ideals the K\"ulshammer ideals.
In joint work \cite{LZZ}
with Liu and Zhou we showed that if the base field is algebraically closed, then
the dimension of the
image of this invariant in the stable centre is also
an invariant under stable equivalences of Morita type.
In joint work \cite{hztame}  with Holm we observed that the K\"ulshammer ideals
behave in a very subtle manner with respect to the deformation parameters.
Using this observation we showed that some of the
parameters are invariants under derived equivalence for certain families
of algebras of dihedral and of semidihedral type. The present paper is
a continuation and completion of \cite{hztame}.

In order to apply the
theory of K\"ulshammer ideals we need to use the symmetrising form
explicitly, and in the present work we progress in avoiding the ad-hoc arguments used in our
previous work to determine the symmetrising form. In this note we
compute the K\"ulshammer ideals for algebras of quaternion type and distinguish
this way the derived equivalence classes of the
algebras with two simple modules. Over algebraically closed fields of
characteristic different from $2$ we can classify completely the derived equivalence classes of
the algebras of quaternion type occurring in Holm's list, except for a case of very small
parameters. If the field is
algebraically closed of characteristic $2$ then we have an almost complete classification
in case of two simple modules. The result in this case is displayed in Corollary~\ref{Finalcorollarytwosimples} and Theorem~\ref{main}.
We also deal with the case of algebras of quaternion type with
three simple modules, where K\"ulshammer ideals distinguish the isomorphism classes of algebras
in characteristic $2$ with parameter $d$ according to whether or not $d$ is a square in $K$.
The invariance of K\"ulshammer ideals under derived or Morita equivalence is proved only
in case $K$ is perfect, which implies that all elements of $K$ are squares when $K$ is of characteristic $2$.
Hence, we cannot say more about this case, and the derived equivalence classification remains open for this class of $20$-dimensional algebras.
Derived equivalent local algebras are actually Morita equivalent
(cf \cite[Proposition 6.7.4]{reptheobuch}), so that the derived equivalence classification of the class of algebras of quaternion type with one simple module coincides with its classification up to isomorphism. Isomorphic algebras have isomorphic K\"ulshammer ideal structure.

For the reader's convenience we give the somewhat technical result for the
class of algebras with two simples here.
Blocks of quaternion type with two simple modules are derived equivalent to an
algebra $A^{k,s}(a,c)$ for parameters $a\in K^\times$ and $c\in K$ and integers
$s\geq 3$ and $k\geq 1$.
\begin{itemize}
\item
In particular, if $K$ is an algebraically closed field of characteristic different from $2$, then
there is $a'\in K^\times$ such that
$A^{k,s}(a,c)\simeq A^{k,s}(a',0)$, and if $(k,s)\neq(1,3)$, then $A^{k,s}(a,c)\simeq A^{k,s}(1,0)$. Moreover, if $A^{k,s}(1,0)$ and $A^{k',s'}(1,0)$ are derived equivalent, then $(k,s)=(k',s')$ or $(k,s)=(s',k')$.
\item
If $K$ is a perfect field of characteristic $2$, we have the following situation.
The algebra $A^{k,s}(a,c)$ is not  derived equivalent to $A^{k,s}(a',0)$ for
any $a,a',c\in K^\times$. If $K$ is algebraically closed of characteristic $2$, and
if $c\neq 0$, then $A^{k,s}(a,c)$ is isomorphic to $A(a'',1)$ for some $a''\in K^\times$, and
if $(k,s)\neq (1,3)$, then $A^{k,s}(a,0)\simeq A^{k,s}(1,0)$.
Further, again for algebraically closed $K$, if
$A^{k,s}(a',c')$ is derived equivalent to $A^{k',s'}(a'',c'')$ then
$(k,s)=(k',s')$ or $(k,s)=(k',s')$.

We do not know for which parameters $a,a'\in K^\times$ we get
$A^{k,s}(a,1)$ is derived equivalent to $A^{k,s}(a',1)$, and we do not know when
$A^{(1,3)}(a,0)$ is derived equivalent to $A^{(1,3)}(a',0)$ for $a,a'\in K^\times$.
\end{itemize}
The K\"ulshammer ideal structure depends in a quite subtle way on the parameters, and
we want to stress the fact that we need to compute the ideals as ideals, and as in \cite{hztame}
it is not sufficient to consider the dimensions only.

The paper is organised as follows. In Section~\ref{Sec-Reynolds} we recall basic facts
about K\"ulshammer ideals and improve the general methods needed to compute the K\"ulshammer
ideal structure for symmetric algebras.
In Section~\ref{quaterniontypesection} we apply the general theory
to algebras of quaternion type,
and we prove our main result Theorem~\ref{main} there.

\subsection*{Acknowledgement}
I wish to thank Oyvind Solberg for giving me during the
Oberwolfach conference ``Hochschild cohomology and applications'' in February 2016
a GAP program to compute the K\"ulshammer ideals. The GAP program \cite{GAP}
uses the package ``qpa'' and encouraged me to study the quaternion type algebras.
I also wish to thank Rachel Taillefer for pointing out the particularity of $(k,s)=(1,3)$
for two simple modules which I forgot to consider in a previous version. I thank the referee for numerous very useful remarks, and in particular for alerting me on some mistake
in the initial version concerning symmetrising forms.

\section{Review on K\"ulshammer ideals and how to compute them}
\label{Sec-Reynolds}

The aim of this section is to briefly give the necessary
background on K\"ulshammer ideals,
as introduced by B. K\"ulshammer
\cite{Ku1}. Morita invariance of K\"ulshammer ideals (then named generalised Reynolds'
ideals) was shown in \cite{BHHKM,HHKM} for perfect fields $K$.
K\"ulshammer ideals were proved to be a derived invariant in \cite{kuelsquest},
were used in \cite{HS,HS2,BSkow} to classify weakly
symmetric algebras of polynomial growth or domestic type up to derived
equivalences, in \cite{hztame} for a derived equivalence classification
of algebras of dihedral or semidihedral type, in \cite{Ln}
for deformed preprojective algebras of type $L$, and in \cite{ST}
for the derived equivalence
classification of certain special biserial algebras. The concept was generalised
to general finite-dimensional algebras in \cite{BHZ}, to an invariant of
Hochschild (co)homology for symmetric algebras \cite{gerstenhaber}
and in \cite{TAHochschild} for general algebras. The image of the K\"ulshammer ideals
in the stable centre were shown to be an
invariant under stable equivalences of Morita type \cite{LZZ,KLZ}.
An overview is given in \cite{Iran,reptheobuch}.

Let $K$ be a field of characteristic $p > 0$.
Any finite-dimensional symmetric $K$-algebra $A$ has
an associative, symmetric, non-degenerate $K$-bilinear form
$\langle -,-\rangle : A
\times A \rightarrow K$. For any $K$-linear subspace $M$ of $A$ we
denote the orthogonal space by $M^{\bot}$ with respect to this form.
Moreover, let $[A,A]$ be the $K$-subspace
of $A$ generated by all commutators $[a,b]:=a b - b a$,
where $a, b \in A$. For any $n \geq 0$ set
$$T_n (A) = \left\{ x \in A \mid x^{p^n} \in [A,A]\right\}.$$
Then, by \cite{Ku1}, for any $n\ge 0$, the orthogonal
space $T_n (A)^{\bot}$ is an ideal of the center
$Z(A)$ of $A$, called {\it $n$-th K\"ulshammer ideal}.
These ideals form a descending sequence
$$Z(A) =  [A,A]^{\perp} = T_0(A)^{\perp} \supseteq T_1(A)^{\perp}
\supseteq T_2(A)^{\perp}
\supseteq \ldots \supseteq T_n(A)^{\perp} \supseteq \ldots$$
with intersection of all ideals $T_n(A)^{\perp}$ for $n\in\N$
being the Reynolds' ideal
$R(A)=Z(A)\cap\soc(A)$.
In \cite{HHKM} it has been shown that if $K$ is perfect, then
the sequence of K\"ulshammer ideals is invariant under Morita equivalences.
Later, it was shown that the sequence of K\"ulshammer ideals
is invariant under derived equivalences, and the image of the
sequence of K\"ulshammer ideals in the stable centre
is invariant under stable equivalences of Morita type.
The following theorem recalls part of what is known.

\begin{Theorem}
\label{prop:zimmermann}
\begin{itemize}\item \cite[Theorem 1]{kuelsquest}
Let $A$ and $B$ be finite-dimensional symmetric algebras over a perfect
field $K$ of positive characteristic $p$.
If $A$ and $B$ are derived equivalent, then
there is an isomorphism $\varphi : Z(A) \rightarrow Z(B)$ between
the centers of $A$ and $B$ such that $\varphi(T_n (A)^{\bot}) = T_n
(B)^{\bot}$ for all positive integers $n$.

\item (cf e.g. \cite[Proposition 6.8.9]{reptheobuch})
Let $A$ and $B$ be derived equivalent finite dimensional $K$-algebras
over a field $K$, which is a splitting field for $A$ and for $B$. Then
the elementary divisors of the Cartan matrices of $A$ and of $B$ coincide. In particular,
the determinant of the Cartan matrices coincides.

\item \cite[Corollary 6.5]{LZZ}
If $A$ and $B$ are stably equivalent of Morita type, and if $K$ is an
algebraically closed field, then $\dim_K(T_n(A)/[A,A])=\dim_K(T_n(B)/[B,B])$.
\end{itemize}
\end{Theorem}

We note that in the proof of \cite[Theorem 1]{kuelsquest} the
hypothesis that $K$ is algebraically closed is never used.
The assumption on the field $K$ to be perfect is sufficient.

The aim of the present note is to show how these derived
invariants can be applied to some subtle questions in the
derived equivalence classifications of algebras of quaternion type.

In order to compute the K\"ulshammer ideals we need a symmetrising form.
However, the K\"uls\-hammer ideals do not depend on the choice of the
symmetrising form if $K$ is perfect (cf \cite[Proof of Claim 3]{kuelsquest}).
We showed in \cite{Ln} (see also \cite{reptheobuch})
that every Frobenius form arises as in the following proposition.

\begin{Prop} \cite{hztame,Ln}
\label{prop:form}
Let $A$ be a basic Frobenius algebra such that $K$ is a splitting field for
$A$, and let $\{e_1,\dots,e_n\}$ be a choice of orthogonal primitive idempotents
with $\sum_{i=1}^n e_i=1$. Then there are bases ${\mathcal B}_{i,j}$ of $e_iAe_j$ such that
${\mathcal B}=\bigcup_{i,j=1}^n{\mathcal B}_{i,j}$ is a basis of $A$ containing
a basis of $\soc(A)$
and such that the following statements hold:
\begin{enumerate}
\item[{(1)}] Define an $K$-linear mapping $\psi$ on the basis elements by
$$
\psi(b)=\left\{
\begin{array}{ll} 1 & \mbox{if $b\in\soc(A)$} \\
                  0 & \mbox{otherwise}
\end{array} \right.
$$
for $b\in {\mathcal B}$.
Then an associative non-degenerate $K$-bilinear
form $\langle-,-\rangle$ for $A$ is given by
$\langle x,y\rangle := \psi(xy).$
\item[{(2)}] Any Frobenius form arises this way for some choice of a basis $\mathcal B$.
\end{enumerate}
\end{Prop}

Note that the hypothesis in \cite{hztame,Ln} is slightly different, however equivalent to the one
given here.

If $A$ is a basic symmetric algebra over an algebraically closed field $K$, then $A=KQ/I$
and we want to determine those bases ${\mathcal B}_s$ of $\soc(A)$ which yield
a symmetric form. This problem is addressed in previous papers
dealing with K\"ulshammer ideals (cf \cite[Remark 2.9]{Ln}, \cite[Remark 3.2]{hztame}).
The following remark indicates a necessary condition for the problem.

\begin{Rem}\label{symmetrisingform}
If $A$ is an indecomposable, basic symmetric algebra over an algebraically closed field $K$ and
let $\{e_1,\dots,e_n\}$ be a choice of orthogonal primitive idempotents
with $\sum_{i=1}^n e_i=1$. Suppose that
${\mathcal B}_s$ is a $K$-basis of $\soc(A)$ and suppose that for each $b\in{\mathcal B}_s$
there is a unique $e_i$ such that $e_ibe_i=b$.
Using \cite[Proposition 2.7.4]{reptheobuch} it is not hard
to see that we can always find a basis
${\mathcal B}_s$ of $\soc(A)$ such that the difference of two elements of ${\mathcal B}_s$
is in the commutator subspace. Moreover, since the elements of ${\mathcal B}_s$
are uniquely determined up to scalars by this property, Proposition~\ref{prop:form} then shows that we can complete the basis ${\mathcal B}_s$ to a basis ${\mathcal B}$ as in the proposition.
If $\psi$ is a $K$-linear map as in
Proposition~\ref{prop:form}, then $\psi([A,A])=0$. In particular, if $b,b'\in \soc(A)$ with
$b-b'\in[A,A]$, then $\psi(b)=\psi(b')$.
\end{Rem}

\section{Algebras of quaternion type}

\label{quaterniontypesection}

\subsection{Two simple modules}

Erdmann gave a classification of algebras which could appear as blocks of
tame representation type. These algebras fall in three classes, the
algebras of dihedral, the algebras of semidihedral and the algebras of
quaternion type. Erdmann's classification was up to Morita equivalence.
Holm \cite[Appendix B]{Holmhabil} gave a classification up to derived
equivalence and obtained for non-local algebras of quaternion type two
families, one containing algebras with two simple modules, one containing algebras
with three simple
modules. The algebras in each family share a common quiver, and the relations depend on a
number of parameters.

The quiver for the algebras with two simples is the following.

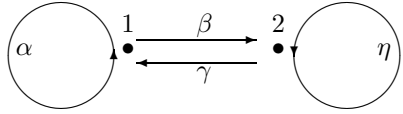
\unitlength1cm
\begin{center}
\begin{picture}(10,2)
\put(2,1){$\bullet$}\put(4,1){$\bullet$}
\put(2,1.3){$1$}\put(4,1.3){$2$}
\put(2.2,1.2){\vector(1,0){1.6}}
\put(3.8,.9){\vector(-1,0){1.6}}
\put(1.2,1){\circle{2}}
\put(1.9,.98){\vector(0,1){.1}}
\put(5,1){\circle{2}}
\put(4.3,1.08){\vector(0,-1){.1}}
\put(.6,1){$\alpha$}
\put(3,1.3){$\beta$}
\put(3,.7){$\gamma$}
\put(5.4,1){$\eta$}
\end{picture}
\end{center}

Let $k\geq 1, s\geq 3$ integers and $a\in K^\times, c\in K$. Then we get an
algebra $Q(2{\mathfrak B})_1^{k,s}(a,c)$ by the above quiver with relations
$$\beta\eta=(\alpha\beta\gamma)^{k-1}\alpha\beta,\;\;\;\;
\eta\gamma=(\gamma\alpha\beta)^{k-1}\gamma\alpha,\;\;\;\;
\alpha^2=a\cdot(\beta\gamma\alpha)^{k-1}\beta\gamma+c\cdot(\beta\gamma\alpha)^k,$$
$$\gamma\beta=\eta^{s-1}, \;\;\;\;\alpha^2\beta=0,\;\;\;\;\gamma\alpha^2=0.$$

\begin{Rem}\label{CentreandCartandet} Using \cite{Holmhabil} we see that
the centre of this algebra is of dimension $k+s+2$ and the Cartan matrix is
$\left(\begin{array}{cc}4k&2k\\2k&k+s\end{array}\right)$ with determinant $4ks$.
Hence, using Theorem~\ref{prop:zimmermann}, if
$D^b(Q(2{\mathfrak B})_1^{k,s}(a,c))\simeq D^b(Q(2{\mathfrak B})_1^{k',s'}(a',c'))$,
then $4ks=4k's'$ and $k+s+2=k'+s'+2$. Therefore $(k+s)^2=(k'+s')^2$ and
$(k-s)^2=(k'-s')^2$, which implies $k=k'$ and $s=s'$, or $k=s'$ and $k'=s$.
\end{Rem}

\begin{Lemma}\label{centreofQ2}
Let $K$ be a field, and let $A^{k,s}(a,c):=Q(2{\mathfrak B})_1^{k,s}(a,c)$.
Then, $Z(A^{k,s}(a,c))$ has a $K$-basis formed by the disjoint union
$$\{\eta-(\alpha\beta\gamma)^{k-1}\alpha\}
\stackrel{\cdot}{\cup}\{\eta^t\;|2\leq t\leq s\}
\stackrel{\cdot}{\cup}\{(\alpha\beta\gamma)^u+(\beta\gamma\alpha)^u+(\gamma\alpha\beta)^u\;|\;
1\leq u\leq k-1\}
\stackrel{\cdot}{\cup}\{1,(\alpha\beta\gamma)^k,\alpha^2\}$$
and is isomorphic, as commutative $K$-algebra, with
$$ K[U,Y,S,T]/(Y^{s+1},U^k-Y^s-2T,S^2,T^2,SY,SU,ST,UY,UT,YT)$$
where
\begin{eqnarray*}
U&:=&(\alpha\beta\gamma)+(\beta\gamma\alpha)+(\gamma\alpha\beta)\\
Y&:=&\eta-(\alpha\beta\gamma)^{k-1}\alpha\\
S&:=&\alpha^2\\
T&:=&(\alpha\beta\gamma)^{k}
\end{eqnarray*}
\end{Lemma}

Proof. First, $\eta-(\alpha\beta\gamma)^{k-1}\alpha$ commutes trivially with $\eta$, since $\eta\alpha=0=\alpha\eta$.
Now
\begin{eqnarray*}
\alpha(\eta-(\alpha\beta\gamma)^{k-1}\alpha)-(\eta-(\alpha\beta\gamma)^{k-1}\alpha)\alpha
&=&\alpha^2(\beta\gamma\alpha)^{k-1}-(\alpha\beta\gamma)^{k-1}\alpha^2\\
&=&(a(\beta\gamma\alpha)^{k-1}\beta\gamma+c(\beta\gamma\alpha)^k)(\beta\gamma\alpha)^{k-1}\\
&&-(\beta\gamma\alpha)^{k-1}(a(\beta\gamma\alpha)^{k-1}\beta\gamma+c(\beta\gamma\alpha)^k)\\
&=&a\left((\beta\gamma\alpha)^{k-1}\beta\gamma(\beta\gamma\alpha)^{k-1}-
(\beta\gamma\alpha)^{2k-2}\beta\gamma\right)\\
&=&a(\beta\gamma\alpha)^{k-1}\left(\beta\gamma(\beta\gamma\alpha)^{k-1}-
(\beta\gamma\alpha)^{k-1}\beta\gamma\right).
\end{eqnarray*}
This is trivially $0$ if $k=1$, and if $k>1$, then $$(\beta\gamma\alpha)^{2k-2}\beta\gamma=\beta(\gamma\alpha\beta)^{2k-2}\gamma
=\beta(\gamma\alpha\beta)^{k}(\gamma\alpha\beta)^{k-2}\gamma=
\beta\eta^s(\gamma\alpha\beta)^{k-2}\gamma=0.$$
Hence
\begin{eqnarray*}
a(\beta\gamma\alpha)^{k-1}\left(\beta\gamma(\beta\gamma\alpha)^{k-1}-
(\beta\gamma\alpha)^{k-1}\beta\gamma\right)
&=&a(\beta\gamma\alpha)^{k-1}\beta\gamma(\beta\gamma\alpha)(\beta\gamma\alpha)^{k-2}\\
&=&a(\beta\gamma\alpha)^{k-1}\beta\eta^{s-1}(\gamma\alpha)(\beta\gamma\alpha)^{k-2}\\
&=&a(\beta\gamma\alpha)^{k-1}\beta\eta^{s-2}(\gamma\alpha\beta)^{k-1}\gamma\alpha^2(\beta\gamma\alpha)^{k-2}\\
&=&0
\end{eqnarray*}
since $\gamma\alpha^2=0$.
The relations
$\beta\eta=(\alpha\beta\gamma)^{k-1}\alpha\beta$ and
$\eta\gamma=(\gamma\alpha\beta)^{k-1}\gamma\alpha=\gamma(\alpha\beta\gamma)^{k-1}\alpha$
show that $\eta-(\alpha\beta\gamma)^{k-1}\alpha$ commutes with $\beta$ and with $\gamma$.
Now, if $k>1$, then $\left(\eta-(\alpha\beta\gamma)^{k-1}\alpha\right)^2=\eta^2$, and if
$k=1$, then $\left(\eta-(\alpha\beta\gamma)^{k-1}\alpha\right)^2=\eta^2+\alpha^2$. Since $\alpha^2\beta=\gamma\alpha^2=0$, it is clear that $\alpha^2$ is central.
Hence $\eta^t$ is central for each $t\geq 2$.
Now,
$$\alpha\beta\gamma\beta=\alpha\beta\eta^{s-1}=\alpha(\alpha\beta\gamma)^{k-1}\alpha\beta\eta^{s-2}=0$$
since $\alpha^2\beta=0$. Likewise  $\gamma\beta\gamma\alpha=0$. Hence
$$
\beta U=\beta\left((\alpha\beta\gamma)+(\beta\gamma\alpha)+(\gamma\alpha\beta)\right)=
\beta\gamma\alpha\beta=\left((\alpha\beta\gamma)+(\beta\gamma\alpha)+(\gamma\alpha\beta)\right)\beta=U\beta
$$
and
$$
\gamma U=\gamma\left((\alpha\beta\gamma)+(\beta\gamma\alpha)+(\gamma\alpha\beta)\right)
=\gamma\alpha\beta\gamma
=\gamma\left((\alpha\beta\gamma)+(\beta\gamma\alpha)+(\gamma\alpha\beta)\right)=
U\gamma.
$$
Now,
\begin{eqnarray*}
\eta U&=&\eta\left((\alpha\beta\gamma)+(\beta\gamma\alpha)+(\gamma\alpha\beta)\right)\\
&=&\eta(\gamma\alpha\beta)=(\gamma\alpha\beta)^{k-1}\gamma\alpha^2\beta\\
&=&0\\
&=&\gamma\alpha(\alpha\beta\gamma)^{k-1}\alpha\beta\\
&=&\gamma\alpha\beta\eta\\
&=&\left((\alpha\beta\gamma)+(\beta\gamma\alpha)+(\gamma\alpha\beta)\right)\eta\\
&=&U\eta
\end{eqnarray*}
and
$$\alpha U=\alpha\left((\alpha\beta\gamma)+(\beta\gamma\alpha)+(\gamma\alpha\beta)\right)
=\alpha\beta\gamma\alpha=\left((\alpha\beta\gamma)+(\beta\gamma\alpha)+(\gamma\alpha\beta)\right)\alpha=
U\alpha.
$$
Hence $U$ is central, and then we only need to compute $U^u$ to get the result. Finally, socle
elements in basic symmetric algebras over splitting fields
are always central, and $1$ is central of course. We know by \cite{Holmhabil}
that the centre is $(2+k+s)$-dimensional,
and obtain therefore the result.
\dickebox

\begin{Rem}  \label{remarkparameters}
Erdmann and Skowro\'nski show in \cite[Lemma 5.7]{ErdSkow}
that if $K$ is an algebraically closed field, then
$Q(2{\mathfrak B})_1^{k,s}(a,c)\simeq Q(2{\mathfrak B})_1^{k,s}(1,c')$ for some $c'\in K$
and if $K$ is of characteristic different from $2$, then
$Q(2{\mathfrak B})_1^{k,s}(a,c)\simeq Q(2{\mathfrak B})_1^{k,s}(a,0)$.
We can examine their computations again to get slightly better results. We assume here $k+s>4$.

Suppose that $K$ admits any $k$-th root, i.e. for all $x\in K$ there is
$y\in K$ with $y^k=x$.
We want to simplify the parameters $a,c$. Replace $\alpha$ by $x_\alpha\alpha$,
$\beta$ by $x_\beta\beta$,
$\gamma$ by $x_\gamma\gamma$ and $\eta$ by $x_\eta\eta$ for non zero scalars
$x_\alpha,x_\beta,x_\gamma,x_\eta$.
Then the relations above are equivalent to
\begin{eqnarray*}
x_\eta\beta\eta&=&x_\alpha^k(x_\beta x_\gamma)^{k-1}(\alpha\beta\gamma)^{k-1}\alpha\beta,\\
x_\eta \eta\gamma&=&x_\alpha^{k}(x_\beta x_\gamma)^{k-1}(\gamma\alpha\beta)^{k-1}\gamma\alpha,\\
x_\alpha^2\alpha^2&=&a\cdot x_\alpha^{k-1}(x_\beta x_\gamma)^{k}
(\beta\gamma\alpha)^{k-1}\beta\gamma+c\cdot x_\alpha^{k}(x_\beta x_\gamma)^{k}(\beta\gamma\alpha)^k,\\
x_\gamma x_\beta\gamma\beta&=&x_\eta^{s-1}\eta^{s-1},\\
\alpha^2\beta&=&0,\\
\gamma\alpha^2&=&0.
\end{eqnarray*}
We first choose $x_\beta$ such that $x_\beta x_\gamma=x_\eta^{s-1}$ to get the system
\begin{eqnarray*}
\beta\eta&=&x_\alpha^kx_\eta^{(k-1)(s-1)-1}(\alpha\beta\gamma)^{k-1}\alpha\beta,\\
\eta\gamma&=&x_\alpha^{k}x_\eta^{(k-1)(s-1)-1}(\gamma\alpha\beta)^{k-1}\gamma\alpha,\\
\alpha^2&=&a\cdot x_\alpha^{k-3}x_\eta^{k(s-1)}
(\beta\gamma\alpha)^{k-1}\beta\gamma+c\cdot x_\alpha^{k-2}x_\eta^{k(s-1)}(\beta\gamma\alpha)^k,\\
\gamma\beta&=&\eta^{s-1},\\
\alpha^2\beta&=&0,\\
\gamma\alpha^2&=&0.
\end{eqnarray*}
Then we put $x_\alpha=x_\eta^{-\frac{(k-1)(s-1)-1}{k}}$ and obtain the system
\begin{eqnarray*}
\beta\eta&=&(\alpha\beta\gamma)^{k-1}\alpha\beta,\\
\eta\gamma&=&(\gamma\alpha\beta)^{k-1}\gamma\alpha,\\
\alpha^2&=&a\cdot x_\eta^{-\frac{(k-1)(s-1)-1}{k}\cdot(k-3)+k(s-1)}
(\beta\gamma\alpha)^{k-1}\beta\gamma+
c\cdot x_\eta^{-\frac{(k-1)(s-1)-1}{k}(k-2)+k(s-1)}(\beta\gamma\alpha)^k,\\
\gamma\beta&=&\eta^{s-1},\\
\alpha^2\beta&=&0,\\
\gamma\alpha^2&=&0.
\end{eqnarray*}
Now, $-\frac{(k-1)(s-1)-1}{k}\cdot(k-3)+k(s-1)=0$ implies $k=3$ and $s=1$, or $k=1$ and $s=3$,
which are excluded parameters, where the first case is already excluded since the algebra is
defined only for $s\geq 3$, and where both cases are excluded by our hypothesis.
This number can be simplified and we therefore define
$u(k,s):=(k-3)+(4k-3)(s-1)>0$ for our admissible parameters.
Hence if there is an element
$y_a\in K$ such that $y_a^{u(k,s)}=a^{-k}$, then we may choose such a
parameter $x_\eta$, such that we may assume $a=1$. This holds in particular if $K$
is algebraically closed.
We obtain this way $A^{k,s}(a,0)\simeq A^{k,s}(1,0)$ if $K$ is sufficiently big, i.e. there is an element $y_a$ satisfying $y_a^{u(k,s)}=a^{-k}$.
Moreover, since $s\geq 3$ and $k\geq 1$ we get
$$-\frac{(k-1)(s-1)-1}{k}(k-2)+k(s-1)= \frac{(k-2)+(3k-2)(s-1)}k\geq \frac{(k-2)+2(3k-2)}k=\frac{7k-6}k>0.$$
Let
$v(k,s):=(k-2)+(3k-2)(s-1)$. If there is $y_c\in K$ such that $y_c^{v(k,s)}=c^{-k}$, then
we can therefore choose $x_\eta$ such that we can assume $c=1$ if $c\neq 0$.
Again, this is trivially true if $K$ is algebraically closed.

As a consequence, combining our computation and the result \cite[Lemma 5.7]{ErdSkow},
if $K$ is algebraically closed of characteristic different from $2$, then
$Q(2{\mathfrak B})_1^{k,s}(a,c)\simeq Q(2{\mathfrak B})_1^{k,s}(1,0)$.
\end{Rem}

\begin{Theorem}\label{main}
Let $A^{k,s}(a,c)$ be the algebra $Q(2{\mathfrak B})_1^{k,s}(a,c)$
over a field $K$ of characteristic $p$.  Let $a,c\in  K\setminus\{ 0\}$.
We get the following cases.
\begin{enumerate}
\item \label{item1} Suppose $p=2$.
\begin{enumerate}
\item If $k=1$, and
\begin{enumerate}
\item if $s$ is even or if $a$ is not a square in $K$, then
$$\dim_K(T_1^\perp(A^{k,s}(a,c)))=\dim_K(T_1^\perp(A^{k,s}(a,0))),$$
\item if $s$ is odd and $a$ is a square in $K$, then
$$\dim_K(T_1^\perp(A^{k,s}(a,c)))=\dim_K(T_1^\perp(A^{k,s}(a,0)))-1.$$
\end{enumerate}
\item If $k>1$ is odd, and
\begin{enumerate}
\item if $s$ is even and if $c$ is a square in $K$, then
$$\dim_K(T_1^\perp(A^{k,s}(a,c)))=\dim_K(T_1^\perp(A^{k,s}(a,0))),$$
\item if $s$ is odd or if $c$ is not a square in $K$, then
$$\dim_K(T_1^\perp(A^{k,s}(a,c)))=\dim_K(T_1^\perp(A^{k,s}(a,0)))-1.$$
\end{enumerate}
\item If $k$ is even, and
\begin{enumerate}
\item if $c$ is a square in $K$, then
$$\dim_K(T_1^\perp(A^{k,s}(a,c)))=\dim_K(T_1^\perp(A^{k,s}(a,0))),$$
\item if $c$ is not a square in $K$, then
$$\dim_K(T_1^\perp(A^{k,s}(a,c)))=\dim_K(T_1^\perp(A^{k,s}(a,0)))-1.$$
\end{enumerate}
\end{enumerate}
\item \label{item2} Suppose $K$ is a perfect field of characteristic $p=2$.
\begin{enumerate}
\item Then $A^{k,s}(a,0)\simeq A^{k,s}(1,0)$ and
\begin{enumerate}
\item if $k$ and $s$ are even, then
$$Z(A^{k,s}(a,0))/T_1^\perp(A^{k,s}(a,0))\simeq K[U,Y,S]/(U^{k/2}-Y^{s/2},S^2,UY,US,YS),$$
\item if $k>1$ or $s$ is odd, then
$$Z(A^{k,s}(a,0))/T_1^\perp(A^{k,s}(a,0))\simeq K[U,Y,S]/(U^{\lceil k/2\rceil},Y^{\lceil s/2\rceil},S^2,UY,US,YS).$$
\end{enumerate}
\item
If $c\neq 0$, then
\begin{enumerate}
\item if $k$ and $s$ are even, then
$$Z(A^{k,s}(a,c))/T_1^\perp(A^{k,s}(a,c))\simeq K[U,Y]/(U^{k/2}-Y^{s/2},UY),$$
\item if $k>1$ or $s$ is odd, then
$$Z(A^{k,s}(a,c))/T_1^\perp(A^{k,s}(a,c))\simeq K[U,Y]/(U^{\lceil (k+1)/2\rceil},Y^{\lceil (s+1)/2\rceil},UY).$$
\end{enumerate}
\item If $k=1$, then
\begin{enumerate}
\item if ($s$ is odd and $c=0$) or ($s$ is even and $c\neq 0$),
$$Z(A^{1,s}(a,c))/T_1^\perp(A^{1,s}(a,c))\simeq K[Y,S]/(Y^{\lceil (s+1)/2\rceil},S^2,YS),$$
\item if ($s$ is odd and $c\neq 0$) or ($s$ is even and $c= 0$),
$$Z(A^{1,s}(a,c))/T_1^\perp(A^{1,s}(a,c))\simeq K[Y]/(Y^{\lceil s/2\rceil}).$$
\end{enumerate}
\end{enumerate}
\item \label{item3} If $p\neq 2$, and if $K$ is algebraically closed, then $A^{k,s}(a,c)\simeq A^{k,s}(a,c)$ for some $a'\in K^\times,$ and if $(k,s)\neq (1,3)$, then $A^{k,s}(a,c)\simeq A^{k,s}(1,0)$.
\item \label{item4} If $p>3$ or $n\geq 2$, then the dimension of the K\"ulshammer ideal
$T_n(A^{k,s}(a,c))^\perp$ does not depend on
the parameters $a,c$.

\end{enumerate}
Suppose that $K$ is algebraically closed.
Then $A^{k,s}(a,0)\simeq A^{k,s}(1,0)$ if $(k,s)\neq (1,3)$ and if $c\neq 0$ then
$A^{k,s}(a,c)\simeq A^{k,s}(a',1)\simeq A^{k,s}(1,c')$ for some $a',c'\in K^\times$.
\end{Theorem}

\begin{Rem}
Consider the case $s=3$ and $k=1$. The computations in Remark~\ref{remarkparameters}
show that if $K$ is algebraically closed, and $c\neq 0$, then for
each $a'$ there is $a$ such that
$A^{1,3}(a',c)\simeq  A^{1,3}(a,1)$.

This case is quite particular which allows
an alternative argument for distinguishing derived equivalence classes.
If $c=0$, then the relations of  $A^{1,3}(a,0)$ are homogeneous. Therefore
the algebra $A^{1,3}(a,0)$ is graded by path lengths with semisimple degree $0$ component.
A theorem of Rouquier \cite[Theorem 6.1]{Rouqgrad} shows that
if $A^{1,3}(a,0)$ is derived equivalent to another algebra $B$, then the induced
stable
equivalence of Morita type induces a grading on $B$. Moreover,
by \cite[Lemma 5.21]{Rouqgrad}
the degree $0$ component of $A^{1,3}(a,0)$ is of finite global dimension if and only if
the degree $0$ component of $B$ is of finite global dimension.
\end{Rem}

\begin{Rem}
The hypothesis that $K$ is algebraically closed is stronger than required. A more precise,
somewhat technical statement is given at the end of Remark~\ref{remarkparameters}.
\end{Rem}

Proof of Theorem~\ref{main}. The isomorphisms
$A^{k,s}(a,0)\simeq A^{k,s}(1,0)$ for $(k,s)\neq (1,3)$ and if $c\neq 0$ then $A^{k,s}(1,c')\simeq A^{k,s}(a,c)\simeq A^{k,s}(a',1)$
for some $a',c'\in K^\times$
follow from Remark~\ref{remarkparameters}.

Define the following subsets of $Q(2{\mathfrak B})_1^{k,s}(a,c)$.
$${\mathcal B}_1:=\{\alpha(\beta\gamma\alpha)^n,(\beta\gamma\alpha)^n\beta\gamma,
(\beta\gamma\alpha)^m, e_1(\alpha\beta\gamma)^\ell\;|\;0\leq\ell\leq k,
1\leq m\leq k-1,0\leq n\leq k-1\},$$
$${\mathcal B}_2:=\{e_2\eta^t,(\gamma\alpha\beta)^m\;|\;0\leq t\leq s,1\leq m\leq k-1\},$$
$${\mathcal B}_3:=\{(\beta\gamma\alpha)^n\beta,\alpha(\beta\gamma\alpha)^n\beta \;|\;
0\leq n\leq k-1\},$$
$${\mathcal B}_4:=\{(\gamma\alpha\beta)^n\gamma,(\gamma\alpha\beta)^n\gamma\alpha \;|\;
0\leq n\leq k-1\}.$$
The (disjoint) union of these sets forms a basis of $Q(2{\mathfrak B})_1^{k,s}(a,c)$, using
the known Cartan matrix of $Q(2{\mathfrak B})_1^{k,s}(a,c)$.

As a next step we need to compute the commutator space. Clearly, non
closed paths are commutators, since if $e_ipe_j\neq 0$ for some path $p$
and $e_i\neq e_j$, then
$p=e_ip-pe_i$. Hence ${\mathcal B}_3\cup {\mathcal B}_4\subseteq [A^{k,s}(a,c),A^{k,s}(a,c)]$.
Moreover, $$\alpha(\beta\gamma\alpha)^n=
\alpha(\beta\gamma\alpha)^n-(\beta\gamma\alpha)^n\alpha\in [A(a,c),A(a,c)]\;\forall n\geq 1,$$
$$(\beta\gamma\alpha)^m=(\alpha\beta\gamma)^m=
(\gamma\alpha\beta)^m\in A^{k,s}(a,c)/[A^{k,s}(a,c),A^{k,s}(a,c)]\;\forall m\geq 0,$$
$$(\beta\gamma\alpha)^m\beta\gamma=\gamma(\beta\gamma\alpha)^m\beta=0
\in A^{k,s}(a,c)/[A^{k,s}(a,c),A^{k,s}(a,c)]\;\forall m\geq 1,$$
and hence $${\mathcal B}_{comm}:=\{\alpha,(\alpha\beta\gamma)^m,
\eta^t,e_1,e_2\;|\;1\leq t\leq s-1,1\leq m\leq k\}$$
is a generating set of $A^{k,s}(a,c)/[A^{k,s}(a,c),A^{k,s}(a,c)]$. Since the dimension of
the centre of $A^{k,s}(a,c)$ equals the dimension of $A^{k,s}(a,c)/[A^{k,s}(a,c),A^{k,s}(a,c)]$,
the algebra $A^{k,s}(a,c)$ being symmetric, and both are of dimension $2+k+s$,
by \cite{Holmhabil}, we get that the classes represented by the elements
${\mathcal B}_{comm}$ form actually a basis
of $A^{k,s}(a,c)/[A^{k,s}(a,c),A^{k,s}(a,c)]$.

We only need to work in $A^{k,s}(a,c)/[A^{k,s}(a,c),A^{k,s}(a,c)]$, and therefore
we need to consider linear combinations
of paths in ${\mathfrak B}_{comm}$ only. We can omit idempotents,
since computing modulo the radical these idempotents
remain idempotents, and are hence never nilpotent modulo commutators.  Hence
we only need to consider linear combinations of elements in the set
$$\{\eta^t,\alpha,(\beta\gamma\alpha)^m\;|\;1\leq t\leq s-1,1\leq m\leq k\}.$$

We deal with the case $p=2$.
Let hence  $p=2$.
In the commutator quotient squaring is semilinear
(cf e.g. \cite{Ku1},\cite[Lemma 2.9.3]{reptheobuch}).

{\bf If $k>1$ is odd,} then
\begin{eqnarray*}
0&=&\left(\sum_{t=1}^{s-1}x_t\eta^t+u\alpha+\sum_{m=1}^ky_m(\beta\gamma\alpha)^m\right)^2\\
&=&\sum_{t=1}^{s-1}x_t^2\eta^{2t}+
u^2c(\beta\gamma\alpha)^k+\sum_{m=1}^ky_m^2(\beta\gamma\alpha)^{2m}\\
&=&\sum_{1\leq t\leq s/2}x_t^2\eta^{2t}+
u^2c(\beta\gamma\alpha)^k+\sum_{m=1}^{(k-1)/2}y_m^2(\beta\gamma\alpha)^{2m},
\end{eqnarray*}
which implies $x_t=0$ for all $t\leq \frac{s}2$, $y_m=0$ for all $m\leq \frac{k-1}2$.
If $c=0$ then there is no other constraint.
Suppose $c\neq 0$.
Then, if $s$ is odd we get $u=0$. If $s$ is even, then
$x_{s/2}^2+cu^2=0$ which has a non trivial solution if and only if $c$ is a square in $K$.
Hence, computing in $A^{k,s}(a,c)/[A^{k,s}(a,c),A^{k,s}(a,c)]$ we get
$$T_1(A^{k,s}(a,c))=\left\{\begin{array}{ll}
\langle \alpha,\eta^t,(\beta\gamma\alpha)^m\;|\;t>\frac s2,m>\frac{k-1}2\rangle_K&\mbox{ if $c=0$}\\
\langle \eta^t,(\beta\gamma\alpha)^m\;|\;t>\frac s2,m>\frac{k-1}2\rangle_K&\mbox{ if $s$ is odd and $c\neq 0$}\\
\langle \eta^t,(\beta\gamma\alpha)^m\;|\;
t>\frac s2,m>\frac{k-1}2\rangle_K&\mbox{ if $s$ is even and $0\neq c\not\in K^2$}\\
\langle \eta^t,(\beta\gamma\alpha)^m,\eta^{s/2}+d\alpha\;|\;
t>\frac s2,m>\frac{k-1}2\rangle_K&\mbox{ if $s$ is even and $0\neq c=d^2$}
\end{array}\right.$$

{\bf If $k$ is even,} then
\begin{eqnarray*}
0&=&\left(\sum_{t=1}^{s-1}x_t\eta^t+u\alpha+\sum_{m=1}^ky_m(\beta\gamma\alpha)^m\right)^2\\
&=&\sum_{t=1}^{s-1}x_t^2\eta^{2t}+
u^2c(\beta\gamma\alpha)^k+\sum_{m=1}^ky_m^2(\beta\gamma\alpha)^{2m}\\
&=&\sum_{1\leq t\leq (s-1)/2}x_t^2\eta^{2t}+
u^2c(\beta\gamma\alpha)^k+\sum_{m=1}^{k/2}y_m^2(\beta\gamma\alpha)^{2m},
\end{eqnarray*}
which implies $y_m=0$ for $1\leq m<k/2$ and $x_t=0$
for $t\leq \frac{s-1}2$.

If $c=0$, then $x_{s/2}+y_{k/2}=0$ in case $s$ is even, and $y_{k/2}=0$ in case $s$ is odd.

Suppose $c\neq 0$.
If $s$ is odd, then $y_{k/2}^2+cu^2=0$, and if $s$ is even, then $y_{k/2}^2+x_{s/2}^2+cu^2=0$.
Again $y_{k/2}^2+cu^2=0$ and
$y_{k/2}^2+x_{s/2}^2+cu^2=0$ has non zero solutions if and only if $c$ is a square.

Computing in $A^{k,s}(a,c)/[A^{k,s}(a,c),A^{k,s}(a,c)]$ we get
$$T_1(A(a,c))=\left\{\begin{array}{ll}
\langle \alpha,\eta^t,(\beta\gamma\alpha)^m\;|\;t>\frac{s-1}2,m>\frac{k}2\rangle_K&
\mbox{ if $s$ is odd and $c=0$}\\
\langle \alpha,\eta^{s/2}+(\beta\gamma\alpha)^{k/2},\eta^t,(\beta\gamma\alpha)^m\;|\;
t>\frac{s}2,m>\frac{k}2\rangle_K&\mbox{ if $s$ is even and $c=0$}\\
\langle \eta^t,(\beta\gamma\alpha)^m\;|\;t>\frac{s-1}2,m>\frac{k}2\rangle_K&
\mbox{ if $s$ is odd and $0\neq c\not\in K^2$}\\
\langle \eta^t,(\beta\gamma\alpha)^m,(\beta\gamma\alpha)^{k/2}+d\alpha\;|\;
t>\frac{s-1}2,m>\frac{k}2\rangle_K&\mbox{ if $s$ is odd and $0\neq c=d^2$}\\
\langle \eta^t,(\beta\gamma\alpha)^m,
(\beta\gamma\alpha)^{k/2}+\eta^{s/2}\;|\;t>\frac{s}2,m>\frac{k}2\rangle_K&
\mbox{ if $s$ is even and $0\neq c\not\in K^2$}\\
\langle \eta^t,(\beta\gamma\alpha)^m,(\beta\gamma\alpha)^{k/2}+\eta^{s/2},\eta^{s/2}+d\alpha\;|\;
t>\frac s2,m>\frac{k}2\rangle_K&\mbox{ if $s$ is even and $0\neq c=d^2$}
\end{array}\right.$$

{\bf If $k=1$,} then, since $\beta\gamma=\gamma\beta=\eta^{s-1}$ in
the commutator quotient,
\begin{eqnarray*}
0=\left(\sum_{t=1}^{s-1}x_t\eta^t+u\alpha+y_1(\beta\gamma\alpha)\right)^2
&=&\sum_{t=1}^{s-1}x_t^2\eta^{2t}+u^2a\eta^{s-1}+
u^2c(\beta\gamma\alpha)
\end{eqnarray*}
which implies $x_t=0$ for $1\leq t\leq \frac{s-2}2$.

If $c=0$, then $x_{s/2}=0$ in case $s$ is even, and $x_{(s-1)/2}^2+au^2=0$ in case $s$ is odd.
This last equation has non zero solutions if and only if $a\in K^2$.

Suppose $c\neq 0$. If $s$ is even, then
$x_{s/2}^2+cu^2=0$. This has non zero solutions if and only if $c\in K^2$.
If $s$ is odd, then $cu^2=0$ and $x_{(s-1)/2}^2+au^2=0$. Hence
$s$ odd implies $u=0=x_{(s-1)/2}$. Computing again in $A^{k,s}(a,c)/[A^{k,s}(a,c),A^{k,s}(a,c)]$,
$$T_1(A^{k,s}(a,c))=\left\{\begin{array}{ll}
\langle \eta^{(s-1)/2}+b\alpha,\eta^t,(\beta\gamma\alpha)\;|\;
t>\frac{s-1}2\rangle_K&\mbox{ if $s$ is odd, $a=b^2$ and $c=0$}\\
\langle \eta^t,(\beta\gamma\alpha)\;|\;t>\frac{s-1}2\rangle_K&\mbox{ if $s$ is odd, $a\not\in K^2$ and $c=0$}\\
\langle \eta^t,(\beta\gamma\alpha)\;|\;t>\frac{s}2\rangle_K&\mbox{ if $s$ is even and $c=0$}\\
\langle \eta^t,(\beta\gamma\alpha)\;|\;t>\frac{s-1}2\rangle_K&\mbox{ if $s$ is odd and $c\neq 0$}\\
\langle \eta^t,(\beta\gamma\alpha)\;|\;t>\frac{s}2\rangle_K&\mbox{ if $s$ is even and $0\neq c\not\in K^2$}\\
\langle \eta^t,(\beta\gamma\alpha),\eta^{s/2}+d\alpha\;|\;t>\frac s2\rangle_K&\mbox{ if $s$ is even and $0\neq c=d^2$}
\end{array}\right.$$

It is easy to see that computing $T_n(A^{k,s}(a,c))/[A^{k,s}(a,c),A^{k,s}(a,c)]$
for $n\geq 2$ (and any $p\geq 2$ in this case)
yields expressions which are independent of $a,c$.

\bigskip

In order to compute the K\"ulshammer ideal $T_1(A^{k,s}(a,c))^\perp$ we need to give the symmetrising form of $A^{k,s}(a,c)$.
Recall that we have a basis ${\mathcal B}=\bigcup_{i=1}^4{\mathcal B}_i$
of $A^{(k,s)}(a,c)$ given by
$${\mathcal B}_1:=\{\alpha(\beta\gamma\alpha)^n,(\beta\gamma\alpha)^n\beta\gamma,
(\beta\gamma\alpha)^m, e_1(\alpha\beta\gamma)^\ell\;|\;0\leq\ell\leq k,
1\leq m\leq k-1,0\leq n\leq k-1\},$$
$${\mathcal B}_2:=\{e_2\eta^t,(\gamma\alpha\beta)^m\;|\;0\leq t\leq s,1\leq m\leq k-1\},$$
$${\mathcal B}_3:=\{(\beta\gamma\alpha)^n\beta,\alpha(\beta\gamma\alpha)^n\beta \;|\;
0\leq n\leq k-1\},$$
$${\mathcal B}_4:=\{(\gamma\alpha\beta)^n\gamma,(\gamma\alpha\beta)^n\gamma\alpha \;|\;
0\leq n\leq k-1\}.$$
We define a trace map $$A(a,c)\stackrel{\psi}{\lra} K$$ by
$$\psi(\eta^s)=\psi((\alpha\beta\gamma)^k)=1,\mbox{ and }\psi(x)=0\mbox{ if }x \mbox{ is a path in the quiver such that } x\in {\mathcal B}\setminus\soc(A^{k,s}(a,c)).$$
Note that $(\alpha\beta\gamma)^k=\beta\eta\gamma=(\beta\gamma\alpha)^k$.
Remark~\ref{symmetrisingform} indicates that $\psi$ should coincide on these
socle elements for $\psi$ to define a symmetric form. Indeed,
$\eta^s=\eta\gamma\beta=(\gamma\alpha\beta)^k$
and hence
$$\eta^s-(\alpha\beta\gamma)^k=(\gamma\alpha\beta)^k-(\alpha\beta\gamma)^k=
[\gamma,\alpha\beta(\gamma\alpha\beta)^{k-1}]$$
is a commutator.
We need to prove that $\psi(c_1c_2)=\psi(c_2c_1)$ for all elements $c_1,c_2\in{\mathcal B}$.

\noindent
{\bf Case $(c_1,c_2)\in {\mathcal B}_1\times {\mathcal B}_1$:}
We obtain
$\alpha(\beta\gamma\alpha)^{n_1}\cdot \alpha(\beta\gamma\alpha)^{n_2}=0$ if $n_1+n_2>0$, and the case $n_1=n_2=0$ is clearly symmetric.
$$\alpha(\beta\gamma\alpha)^{n_1}\cdot (\beta\gamma\alpha)^{n_2}\beta\gamma=
(\alpha\beta\gamma)^{n_1+n_2+1}=
(\beta\gamma\alpha)^{n_2}\beta\gamma\cdot \alpha(\beta\gamma\alpha)^{n_1},$$
$$\alpha(\beta\gamma\alpha)^{n_1}\cdot (\beta\gamma\alpha)^{m_2}=
(\alpha\beta\gamma)^{n_1+m_2}\alpha\in({\mathcal B}_1\cup\{0\})\setminus\soc(A^{(k,s)}(a,c)$$
is mapped to $0$ by $\psi$, and $m_2>0$ implies
$$(\beta\gamma\alpha)^{m_2}\cdot \alpha(\beta\gamma\alpha)^{n_1}=0.$$
Now, if $\ell_2>0$, then
$$\alpha(\beta\gamma\alpha)^{n_1}\cdot (\alpha\beta\gamma)^{\ell_2}=0$$
and
$$(\alpha\beta\gamma)^{\ell_2}\cdot\alpha(\beta\gamma\alpha)^{n_1}=
(\alpha\beta\gamma)^{\ell_2+n_1}\alpha\in({\mathcal B}_1\cup\{0\})\setminus\soc(A^{(k,s)}(a,c)$$
is mapped to $0$ by $\psi$. If $\ell_2=0$, then the two elements commute trivially.
$$(\beta\gamma\alpha)^{n_1}\beta\gamma\cdot (\beta\gamma\alpha)^{n_2}\beta\gamma=0=
(\beta\gamma\alpha)^{n_2}\beta\gamma\cdot (\beta\gamma\alpha)^{n_1}\beta\gamma$$
and
$$(\beta\gamma\alpha)^{n_1}\beta\gamma\cdot (\beta\gamma\alpha)^{m_2}=0$$
whereas
$$(\beta\gamma\alpha)^{m_2}\cdot(\beta\gamma\alpha)^{n_1}\beta\gamma=
(\beta\gamma\alpha)^{m_2+n_1}\beta\gamma\in({\mathcal B}_1\cup\{0\})\setminus\soc(A^{(k,s)}(a,c)$$
is mapped to $0$ by $\psi.$ If $\ell_2>0$, then
$$(\beta\gamma\alpha)^{n_1}\beta\gamma\cdot (\alpha\beta\gamma)^{\ell_2}=
(\beta\gamma\alpha)^{n_1+\ell_2}\beta\gamma\in({\mathcal B}_1\cup\{0\})\setminus\soc(A^{(k,s)}(a,c)$$
is mapped to $0$ by $\psi$, whereas
$$(\alpha\beta\gamma)^{\ell_2}\cdot(\beta\gamma\alpha)^{n_1}\beta\gamma=0.$$
Clearly $e_1$ commutes with $(\beta\gamma\alpha)^{n_1}\beta\gamma$.
Now, trivially
$$(\beta\gamma\alpha)^{m_1}\cdot (\beta\gamma\alpha)^{m_2}=(\beta\gamma\alpha)^{m_2}\cdot (\beta\gamma\alpha)^{m_1}$$
and
$$(\alpha\beta\gamma)^{\ell_1}\cdot (\alpha\beta\gamma)^{\ell_2}=(\alpha\beta\gamma)^{\ell_2}\cdot (\alpha\beta\gamma)^{\ell_1}.$$
Finally, if $\ell_1>0$ then
$$(\alpha\beta\gamma)^{\ell_1}\cdot (\beta\gamma\alpha)^{m_2}=0
= (\beta\gamma\alpha)^{m_2}\cdot(\alpha\beta\gamma)^{\ell_1}.$$

\noindent
{\bf Case $(c_1,c_2)\in {\mathcal B}_1\times {\mathcal B}_2$:}
Since ${\mathcal B}_1\subseteq e_1A^{k,s}(a,c)e_1$, and since
${\mathcal B}_2\subseteq e_2A^{k,s}(a,c)e_2$ we get $\psi(c_1c_2)=\psi(c_2c_1)=0$ for
$c_1\in{\mathcal B_1}$ and $c_2\in{\mathcal B_2}$.

\noindent
{\bf Case $(c_1,c_2)\in {\mathcal B}_1\times {\mathcal B}_3$:}
Since ${\mathcal B}_1\subseteq e_1A^{k,s}(a,c)e_1$, and since
${\mathcal B}_3\subseteq e_1A^{k,s}(a,c)e_2$ we get $c_1c_2=0$ and $c_2c_1\in e_1A^{k,s}(a,c)e_3$ for
$c_1\in{\mathcal B_3}$ and $c_2\in{\mathcal B_1}$. Non closed paths are mapped to $0$ by $\psi$.

\noindent
{\bf Case $(c_1,c_2)\in {\mathcal B}_1\times {\mathcal B}_4$:}
Since ${\mathcal B}_4\subseteq e_2A^{k,s}(a,c)e_1$ the same arguments as in the case $(c_1,c_2)\in {\mathcal B}_1\times {\mathcal B}_3$ apply.

\noindent
{\bf Case $(c_1,c_2)\in {\mathcal B}_2\times {\mathcal B}_2$:}
Clearly $\eta^{t_1}\cdot\eta^{t_2}=\eta^{t_2}\cdot\eta^{t_1}$ and
$(\gamma\alpha\beta)^{m_1}\cdot(\gamma\alpha\beta)^{m_2}=
(\gamma\alpha\beta)^{m_2}\cdot(\gamma\alpha\beta)^{m_1}.$
Moreover, if $t>0$ then
$$\eta^t\cdot(\gamma\alpha\beta)^m=0=(\gamma\alpha\beta)^m\cdot\eta^t.$$
If $t=0$, then trivially $\eta^t\cdot(\gamma\alpha\beta)^m=(\gamma\alpha\beta)^m\cdot\eta^t.$

\noindent
{\bf Case $(c_1,c_2)\in {\mathcal B}_2\times {\mathcal B}_3$:}
Since then $c_1c_2$ and $c_2c_1$ are non closed paths, the same arguments as in the case
$(c_1,c_2)\in {\mathcal B}_1\times {\mathcal B}_3$ apply.

\noindent
{\bf Case $(c_1,c_2)\in {\mathcal B}_2\times {\mathcal B}_4$:}
Again since then $c_1c_2$ and $c_2c_1$ are non closed paths, the same arguments as in the case
$(c_1,c_2)\in {\mathcal B}_1\times {\mathcal B}_3$ apply.

\noindent
{\bf Case $(c_1,c_2)\in {\mathcal B}_3\times {\mathcal B}_3$:}
Then $c_1c_2=0=c_2c_1$.

\noindent
{\bf Case $(c_1,c_2)\in {\mathcal B}_3\times {\mathcal B}_4$:}
$$(\beta\gamma\alpha)^{n_1}\beta\cdot (\gamma\alpha\beta)^{n_2}\gamma=
(\beta\gamma\alpha)^{n_1+n_2}\beta\gamma\in({\mathcal B}_1\cup\{0\})\setminus\soc(A^{(k,s)}(a,c)$$
is mapped to $0$ by $\psi$,
and
$$(\gamma\alpha\beta)^{n_2}\gamma\cdot (\beta\gamma\alpha)^{n_1}\beta=
(\gamma\alpha\beta)^{n_2}(\gamma\beta)(\gamma\alpha\beta)^{n_1}=0.$$
Now,
$$(\beta\gamma\alpha)^{n_1}\beta\cdot (\gamma\alpha\beta)^{n_2}\gamma\alpha-
(\gamma\alpha\beta)^{n_2}\gamma\alpha\cdot (\beta\gamma\alpha)^{n_1}\beta=
(\beta\gamma\alpha)^{n_1+n_2+1}-(\gamma\alpha\beta)^{n_1+n_2+1}.$$
and the value of $\psi$ on each of the summands is equal.
$$\alpha(\beta\gamma\alpha)^{n_1}\beta\cdot(\gamma\alpha\beta)^{n_2}\gamma=
(\alpha\beta\gamma)^{n_1+n_2+1}$$
and
$$(\gamma\alpha\beta)^{n_2}\gamma\cdot\alpha(\beta\gamma\alpha)^{n_1}\beta=
(\gamma\alpha\beta)^{n_2+n_1+1}$$
both have identical values under $\psi.$
Finally
$$\alpha(\beta\gamma\alpha)^{n_1}\beta\cdot(\gamma\alpha\beta)^{n_2}\gamma\alpha=
(\alpha\beta\gamma)^{n_1+n_2+1}\alpha\in({\mathcal B}_1\cup\{0\})\setminus\soc(A^{(k,s)}(a,c)
$$
is mapped to $0$ by $\psi$
and
$$(\gamma\alpha\beta)^{n_2}\gamma\alpha\cdot\alpha(\beta\gamma\alpha)^{n_1}\beta
=0.$$

\noindent
{\bf Case $(c_1,c_2)\in {\mathcal B}_4\times {\mathcal B}_4$:}
Then $c_1c_2=0=c_2c_1$.

Altogether this shows that $\psi$ is symmetric. The fact that $\psi$ defines a non degenerate
form follows as in the proof of Proposition~\ref{prop:form}. For the reader's convenience
we recall the short argument. Suppose that the form defined by $\psi$ is degenerate. Then
there is a $0\neq x\in A^{k,s}(a,c)$ with $\psi(xy)=0$ for all $y$, and since $1=e_1+e_2$ there is a
primitive idempotent $e\in\{e_1,e_2\}$ of $A^{k,s}(a,c)$ such that we may suppose that
$x\in eA^{k,s}(a,c)$. Let $S$ be a simple submodule of
$xA^{k,s}(a,c)$ and there is $y$ such that
$0\neq s=xy\in S$. Since $S\leq eA^{k,s}(a,c)$ is one-dimensional,
and included in the socle, and since $e{\mathcal B}e$ contains a basis of $S$ we get $\psi(xy)\neq 0$.
The form defined by $\psi$ is trivially associative. Hence $\psi$ defines
a symmetrising form.

\medskip

We come to the main body of the proof.
Recall from Lemma~\ref{centreofQ2} that $\dim_K(Z(A^{k,s}(a,c)))=k+s+2$.
We proceed case by case.

\medskip

{\bf $k>1$ odd and $c=0$}: Recall that in this case
$$T_1(A^{k,s}(a,c))/[A^{k,s}(a,c),A^{k,s}(a,c)]=
\left\langle\alpha,\eta^t,(\beta\gamma\alpha)^m\;|\;t>\frac s2; m>\frac{k-1}2
\right\rangle_K.$$
Hence
\begin{eqnarray*}\dim_K(T_1(A^{k,s}(a,0))/[A^{k,s}(a,0),A^{k,s}(a,0)])&=&
\left\{\begin{array}{ll}1+\frac s2+\frac{k+1}2-1&\mbox{ if $s$ is even}\\
1+\frac{s+1}2+\frac{k+1}2-1&\mbox{ if $s$ is odd}
\end{array}\right.\\
&=&\left\{\begin{array}{ll}\frac s2+\frac{k+1}2&\mbox{ if $s$ is even}\\
\frac{s+1}2+\frac{k+1}2&\mbox{ if $s$ is odd}
\end{array}\right.
\end{eqnarray*}
observing that $(\beta\gamma\alpha)^k-\eta^s\in [A^{k,s}(a,c),A^{k,s}(a,c)]$.
Therefore
\begin{eqnarray*}
\dim_K(T_1(A^{k,s}(a,0))^\perp)&=&k+s+2-\left\{\begin{array}{ll}\frac s2+\frac{k+1}2&\mbox{ if $s$ is even}\\
\frac{s+1}2+\frac{k+1}2&\mbox{ if $s$ is odd}
\end{array}\right.\\
&=&\left\{\begin{array}{ll}\frac s2+1+\frac{k+1}2&\mbox{ if $s$ is even}\\
\frac{s+1}2+\frac{k+1}2&\mbox{ if $s$ is odd}
\end{array}\right.
\end{eqnarray*}

But, in case $s$ is even,
$$\left\{\eta^t,
(\alpha\beta\gamma)^u+(\beta\gamma\alpha)^u+(\gamma\alpha\beta)^u,(\alpha\beta\gamma)^k\;|\;
u\geq\frac{k+1}2, t\geq \frac s2\right\}\subseteq T_1(A(a,c))^\perp,$$
and in case $s$ is odd,
$$\left\{\eta^t,
(\alpha\beta\gamma)^u+(\beta\gamma\alpha)^u+(\gamma\alpha\beta)^u,(\alpha\beta\gamma)^k\;|\;
u\geq\frac{k+1}2, t\geq \frac{s+1}2\right\}\subseteq T_1(A(a,c))^\perp.$$
This is a basis of a subspace of the centre of the dimension as required, and hence
the set above is a basis of $T_1(A^{k,s}(a,c))^\perp$. Hence, with these parameters, if $s$ is even then
$$Z(A^{k,s}(a,c))/T_1^\perp(A^{k,s}(a,c))\simeq K[U,Y,S]/(Y^{s/2},U^{(k+1)/2},S^2,YS,US,UY),$$
and if $s$ is odd, then
$$Z(A^{k,s}(a,c))/T_1^\perp(A^{k,s}(a,c))\simeq K[U,Y,S]/(Y^{(s+1)/2},U^{(k+1)/2},S^2,YS,US,UY).$$

{\bf $k>1$ odd, $c\neq 0$, and $s$ is odd}:
Recall
$$T_1(A^{k,s}(a,c))/[A^{k,s}(a,c),A^{k,s}(a,c)]=
\left\langle\eta^t,(\beta\gamma\alpha)^m\;|\;t>\frac s2; m>\frac{k-1}2
\right\rangle_K.$$
In this case
we get $\alpha^2\in T_1(A^{k,s}(a,c))^\perp$, and using the preceding discussion we get that
$$\left\{\alpha^2, (\beta\gamma\alpha)^k,\eta^t,
(\alpha\beta\gamma)^u+(\beta\gamma\alpha)^u+(\gamma\alpha\beta)^u\;|\;
u\geq\frac{k+1}2,t\geq \frac{s+1}2\right\}$$
is a $K$-basis of $T_1(A^{k,s}(a,c))^\perp$.
Hence in this case
$$Z(A^{k,s}(a,c))/T_1^\perp(A^{k,s}(a,c))\simeq K[U,Y]/(Y^{(s+1)/2},U^{(k+1)/2},UY).$$

{\bf $k>1$ odd, $d^2=c\neq 0$, and $s$ is even.}
Recall
$$T_1(A^{k,s}(a,c))/[A^{k,s}(a,c),A^{k,s}(a,c)]=
\left\langle d\alpha+\eta^{s/2},\eta^t,(\beta\gamma\alpha)^m\;|\;t>\frac s2; m>\frac{k-1}2
\right\rangle_K.$$
Then
$$\left(\eta^{s/2}+d\alpha\right)\cdot\left(\eta^{s/2}+\frac d{ca}\alpha^2\right)=
\eta^s+\frac{d^2}{ca}\alpha^3=
\eta^s+\frac c{ca}\alpha^3=\eta^s+\frac 1a\cdot a(\alpha\beta\gamma)^k
$$
and this is mapped to $0$ by $\psi$. Hence,
$$\left\{\frac{d}{ca}\alpha^2+\eta^{s/2}, (\alpha\beta\gamma)^k,
\eta^t, (\alpha\beta\gamma)^u+(\beta\gamma\alpha)^u+(\gamma\alpha\beta)^u\;|\;
u\geq\frac{k-1}2,t\geq \frac s2+1\right\}$$ is a $K$-basis of $T_1(A^{k,s}(a,c))^\perp$.
Therefore in this case
$$Z(A^{k,s}(a,c))/T_1^\perp(A^{k,s}(a,c))\simeq K[U,Y]/(Y^{\frac s2+1},U^{(k+1)/2},UY).$$

{\bf $k$ even and $c=0$ and $s$ is odd}.
Recall
$$T_1(A^{k,s}(a,c))/[A^{k,s}(a,c),A^{k,s}(a,c)]=
\left\langle d\alpha+(\beta\gamma\alpha)^{k/2},\eta^t,(\beta\gamma\alpha)^m\;|\;t>\frac {s-1}2;
m>\frac{k}2
\right\rangle_K.$$
Then the discussion of the case $k>1$ odd and $c=0$ shows that
$$\left\{\eta^t,(\alpha\beta\gamma)^k,
(\alpha\beta\gamma)^u+(\beta\gamma\alpha)^u+(\gamma\alpha\beta)^u\;|\;
u\geq\frac{k}2,t\geq \frac {(s+1)}2\right\}$$ is a $K$-basis of $T_1(A^{k,s}(a,c))^\perp.$
Hence in this case
$$Z(A^{k,s}(a,c))/T_1^\perp(A^{k,s}(a,c))\simeq K[U,Y,S]/(Y^{ (s+1)/2},U^{k/2},S^2,YS,US,UY).$$

{\bf $k$ even and $c=0$ and $s$ is even}.
Recall
$$T_1(A^{k,s}(a,c))/[A^{k,s}(a,c),A^{k,s}(a,c)]=
\left\langle \alpha,\eta^{s/2}+(\beta\gamma\alpha)^{k/2},\eta^t,(\beta\gamma\alpha)^m\;|\;
t>\frac s2; m>\frac{k}2
\right\rangle_K.$$
Then
$$\left\{(\alpha\beta\gamma)^k,\eta^t,
(\alpha\beta\gamma)^u+(\beta\gamma\alpha)^u+(\gamma\alpha\beta)^u,\eta^{s/2}+(\beta\gamma\alpha)^{k/2}\;|\;
u\geq\frac{k}2+1,t\geq \frac s2+1\right\}$$ is a $K$-basis of $T_1(A^{k,s}(a,c))^\perp.$
Hence in this case
$$Z(A^{k,s}(a,c))/T_1^\perp(A^{k,s}(a,c))\simeq K[U,Y,S]/(Y^{s/2}-U^{k/2},S^2,YS,US,UY).$$

{\bf $k$ even, $0\neq c=d^2$, $s$ odd:}
Recall
$$T_1(A^{k,s}(a,c))/[A^{k,s}(a,c),A^{k,s}(a,c)]=
\left\langle d\alpha+(\beta\gamma\alpha)^{k/2},\eta^t,(\beta\gamma\alpha)^m\;|\;
t>\frac {s-1}2; m>\frac{k}2
\right\rangle_K.$$
Since $\dim(Z(A^{k,s}(a,c)))=k+s+2$, and since
$$\dim(T_1(A^{k,s}(a,c)))/[A^{k,s}(a,c),A^{k,s}(a,c)]=\frac k2+\frac{s+1}2,$$
we get $$\dim(T_1(A^{k,s}(a,c))^\perp)=2+k+s-\frac k2-\frac{s+1}2=\frac k2+\frac{s+1}2+1.$$
Moreover,
$$
\left(\frac{d}{c}\alpha^2+U^{k/2}\right)\cdot\left(d\alpha+(\beta\gamma\alpha)^{k/2}\right)=
\frac{d^2}c\alpha^3+dU^{k/2}\alpha+\frac{d}c(\beta\gamma\alpha)^{k/2}+(\beta\gamma\alpha)^k
$$
and this maps to $0$ by the map $\psi$.
Therefore
$$\left\{(\beta\gamma\alpha)^k,
\frac{d}{c}\alpha^2+(\alpha\beta\gamma)^{k/2}+(\beta\gamma\alpha)^{k/2}+(\gamma\alpha\beta)^{k/2},
\eta^t,(\alpha\beta\gamma)^u+(\beta\gamma\alpha)^u+(\gamma\alpha\beta)^u\;|\;
u\geq\frac{k}2+1,t\geq \frac {s+1}2\right\}$$ is a $K$-basis of $T_1(A^{k,s}(a,c))^\perp,$
and therefore
$$Z(A^{k,s}(a,c))/T_1(A^{k,s}(a,c))^\perp\simeq K[U,Y]/(U^{k/2},Y^{(s+1)/2},UY).$$

{\bf $k$ even, $c=d^2\neq 0$, and $s$ even:}
Recall
$$T_1(A^{k,s}(a,c))/[A^{k,s}(a,c),A^{k,s}(a,c)]=
\left\langle
d\alpha+\eta^{s/2},(\beta\gamma\alpha)^{k/2}+\eta^{s/2},\eta^t,(\beta\gamma\alpha)^m\;|\;
t>\frac s2; m>\frac{k}2
\right\rangle_K.$$
Then
$$\left\{\alpha^2+d\eta^{s/2},(\beta\gamma\alpha)^k,\eta^t,
(\alpha\beta\gamma)^u+(\beta\gamma\alpha)^u+(\gamma\alpha\beta)^u,\eta^{s/2}+(\beta\gamma\alpha)^{k/2}\;|\;
u\geq\frac{k}2+1,t\geq \frac {s}2+1\right\}$$ is a $K$-basis of $T_1(A^{k,s}(a,c))^\perp.$
Hence in this case
$$Z(A^{k,s}(a,c))/T_1^\perp(A^{k,s}(a,c))\simeq K[U,Y]/(Y^{ s/2}-U^{k/2},UY).$$

{\bf If $k=1$ and $c=0$ and $s$ odd:}
Since $\dim(Z(A^{1,s}(a,c))=3+s$, and since $$\dim(T_1(A^{1,s}(a,c))/[A^{1,s}(a,c),A^{1,s}(a,c)])=3+\frac{s-1}2,$$ we obtain $\dim(T_1(A^{1,s}(a,c))^\perp)=\frac{s+1}2$.
Observe that $\eta^s=(\alpha\beta\gamma)+(\beta\gamma\alpha)+(\gamma\alpha\beta)=U$.
Then we get
$$\left\{\beta\gamma\alpha,\eta^t,\;|\;
t\geq \frac {s+1}2\right\}$$ is a $K$-basis of $T_1(A^{1,s}(a,c))^\perp$.
Therefore
$$Z(A^{1,s}(a,c))/T_1^\perp(A^{1,s}(a,c))\simeq K[Y,S]/(Y^{(s+1)/2},S^2,YS).$$

{\bf If $k=1$ and $c=0$ and $s$ even:}
Since $\dim(Z(A^{1,s}(a,c))=3+s$, and since $$\dim(T_1(A^{1,s}(a,c))/[A^{1,s}(a,c),A^{1,s}(a,c)])=1+\frac{s}2,$$ we obtain $\dim(T_1(A^{1,s}(a,c))^\perp)=2+\frac{s}2$.
Hence
$$\left\{\alpha^2,\beta\gamma\alpha,\eta^t,\;|\;
t\geq \frac {s}2\right\}$$ is a $K$-basis of $T_1(A^{1,s}(a,c))^\perp$ and
$$Z(A^{1,s}(a,c))/T_1^\perp(A^{1,s}(a,c))\simeq K[Y]/Y^{s/2}.$$

{\bf If $k=1$ and $c\neq 0$ and $s$ odd:}
Since $\dim(Z(A^{1,s}(a,c))=3+s$, and since $$\dim(T_1(A^{1,s}(a,c))/[A^{1,s}(a,c),A^{1,s}(a,c)])=2+\frac{s-1}2,$$ we obtain $\dim(T_1(A^{1,s}(a,c))^\perp)=1+\frac{s+1}2$.
Hence
$$\left\{\alpha^2,\beta\gamma\alpha,\eta^t,\;|\;
t\geq \frac {s+1}2\right\}$$ is a $K$-basis of $T_1(A^{1,s}(a,c))^\perp$ and
$$Z(A^{1,s}(a,c))/T_1^\perp(A^{1,s}(a,c))\simeq K[Y]/Y^{(s+1)/2}.$$

{\bf If $k=1$ and $c\neq 0$ and $s$ even:}
Since $\dim(Z(A^{1,s}(a,c))=3+s$, and since $$\dim(T_1(A^{1,s}(a,c))/[A^{1,s}(a,c),A^{1,s}(a,c)])=2+\frac{s}2,$$ we obtain $\dim(T_1(A^{1,s}(a,c))^\perp)=1+\frac{s}2$.
Hence
$$\left\{\beta\gamma\alpha,\eta^t,\;|\;
t\geq 1+\frac {s}2\right\}$$ is a $K$-basis of $T_1(A^{1,s}(a,c))^\perp$ and
$$Z(A^{1,s}(a,c))/T_1^\perp(A^{1,s}(a,c))\simeq K[Y,S]/(Y^{(s+2)/2},S^2,YS).$$

\dickebox

In order to be more concise we summarise the results from Theorem~\ref{main} and Remark~\ref{CentreandCartandet}
in case $K$ is algebraically closed in the following corollary.

\begin{Cor}\label{Finalcorollarytwosimples}
Let $K$ be an algebraically closed field of characteristic $p\in\N\cup\{\infty\}$, and let
$a,a',c$ be non-zero elements in $K$, and let $c',c''\in K$.
\begin{itemize}
\item
If $p\neq 2$, then there is $a'\in K^\times$ such that
$Q(2{\mathfrak B})_1^{k,s}(a,c)\simeq Q(2{\mathfrak B})_1^{k,s}(a',0)$,
and if $(k,s)\neq (1,3)$, then
$Q(2{\mathfrak B})_1^{k,s}(a,c)\simeq Q(2{\mathfrak B})_1^{k,s}(1,0)$.
\item
If $p=2$, then  $D^b(Q(2{\mathfrak B})_1^{k,s}(a,c))\not\simeq D^b(Q(2{\mathfrak B})_1^{k,s}(a',0))$.
Moreover, there is $a''\in K^\times$ such that
$Q(2{\mathfrak B})_1^{k,s}(a,c)\simeq Q(2{\mathfrak B})_1^{k,s}(a'',1)$ and if $(k,s)\neq (1,3)$, then $Q(2{\mathfrak B})_1^{k,s}(a',0)\simeq Q(2{\mathfrak B})_1^{k,s}(1,0)$.
\item
For any characteristic of $K$ we get
$$\left(D^b(Q(2{\mathfrak B})_1^{k,s}(a,c''))\simeq D^b(Q(2{\mathfrak B})_1^{k',s'}(a',c'))\right)\Rightarrow
\left(\mbox{$(k=k'$ and $s=s'$) or ($k=s'$ and $s=k')$.}\right)$$
\end{itemize}
\end{Cor}

Proof. The first statement is an immediate consequence of Theorem~\ref{main} item (\ref{item3}) and
\cite[Lemma 5.7]{ErdSkow}.
The second statement follows from Theorem~\ref{main} item (\ref{item2})(a), (\ref{item2})(b), (\ref{item2})(c), and Theorem~\ref{prop:zimmermann}.
Indeed, the isomorphism type of the centre modulo the first K\"ulshammer ideal
differs in case $c=0$ and $c\neq 0$.
More precisely,
the commutative algebras from case (\ref{item2})(a) (i) and (\ref{item2})(a) (ii) are non isomorphic since
the dimensions of the socles of these algebras differ by $1$. Likewise,
the commutative algebras from case (\ref{item2})(b) (i) and (\ref{item2})(b) (ii) are non isomorphic since
the dimensions of the socles of these algebras differ by $1$.
The dimension of the socle of the centre modulo the K\"ulshammer ideal distinguish the algebras also in case (\ref{item2})(c), i.e. $k=1$.
The third statement follows from Remark~\ref{CentreandCartandet}.
\dickebox

\begin{Rem} Let $K$ be an algebraically closed field of characteristic $2$.
We do not know for which pair of parameters $a,a'\in K^\times$ we get that
$Q(2{\mathfrak B})_1^{k,s}(a',1)$ is derived equivalent to $Q(2{\mathfrak B})_1^{k,s}(a,1)$.
We do not know for which parameters $k,s$ the algebras  $Q(2{\mathfrak B})_1^{k,s}(a,c)$
and $Q(2{\mathfrak B})_1^{s,k}(a,c)$ are derived equivalent.
\end{Rem}

\begin{Rem}
The case $p=3$ is special if $K$ is not perfect.
Let $p=3$ and use the notations used in the proof of Theorem~\ref{main}. Then $\alpha^3=a(\beta\gamma\alpha)^m$. In
the commutator quotient taking third power is semilinear
(cf e.g. \cite{Ku1},\cite[Lemma 2.9.3]{reptheobuch}), and therefore
\begin{eqnarray*}
0&=&\left(\sum_{t=1}^{s-1}x_t\eta^t+u\alpha+\sum_{m=1}^ky_m(\beta\gamma\alpha)^m\right)^3\\
&=&\sum_{t=1}^{s-1}x_t^3\eta^{3t}+
u^3a(\beta\gamma\alpha)^k+\sum_{m=1}^ky_m^3(\beta\gamma\alpha)^{3m}\\
&=&\sum_{1\leq t\leq (s-1)/3}x_t^3\eta^{3t}+
u^3a(\beta\gamma\alpha)^k+\sum_{1\leq m\leq k/3}y_m^3(\beta\gamma\alpha)^{3m}.
\end{eqnarray*}
We have again various cases. If $3$ does not divide $k$ and $3$ does not divide $s$, then
$x_t=0$ for all $t\leq s/3$ and $y_m=0$ for all $m\leq k/3$ and $u=0$.
If $3$ does not divide $k$ but $3|s$, then $x_t=0$ for all $t< s/3$ and
$y_m=0$ for all $m\leq k/3$ and $x_{s/3}^3+au^3=0$, which has a non zero
solution if and only if $a$ is a cube. If $3$ divides $k$ and $3$ does not divide $s$, then
$x_t=0$ for all $t\leq s/3$ and $y_m=0$ for all $m< k/3$ and $y_{m/3}^3+au^3=0$,
which has a non zero solution if and only if $a$ is a cube.
If $3$ divides $k$ and $3$ divides $s$, then
$x_t=0$ for all $t< s/3$ and $y_m=0$ for all $m< k/3$ and
$x_{s/3}^3+y_{m/3}^3+au^3=0$, which has a non zero solution if and only if $a$ is a cube.

As seen above, the first K\"ulshammer ideal
detects if the parameter $a$ is a third power in case $k$ or $s$ is divisible by $3$.
This shows that
the isomorphism $A^{k,s}(a,0)\simeq A^{k,s}(1,0)$, which we proved for
algebraically closed base fields,
is false if the base field is not perfect.
\end{Rem}

\subsection{Three simple modules}

Holm shows that there are two families of algebras $Q(3{\mathcal K})^{a,b,c}$ and
$Q(3{\mathcal A})^{2,2}_1(d)$ with three simple modules such that any block with
quaternion defect group and three simple modules is derived equivalent to an algebra
in one of these families. According to \cite{Holmhabil}
the derived classification of the case
$Q(3{\mathcal K})^{a,b,c}$ is complete, whereas the classification for the case
$Q(3{\mathcal A})^{2,2}_1(d)$ is complete up to the scalar $d\in K\setminus\{0,1\}$.

The quiver $3{\mathcal A}$ is

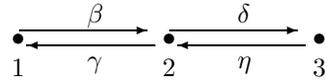
\unitlength1cm
\begin{center}
\begin{picture}(10,2)
\put(0,.9){$\bullet$}\put(2,.9){$\bullet$}\put(4,.9){$\bullet$}
\put(0,.5){$1$}\put(2,.5){$2$}\put(4,.5){$3$}
\put(.1,1.1){\vector(1,0){1.7}}
\put(2.1,1.1){\vector(1,0){1.7}}
\put(1.9,.9){\vector(-1,0){1.7}}
\put(3.9,.9){\vector(-1,0){1.7}}
\put(1,1.2){$\beta$}
\put(3,1.2){$\delta$}
\put(1,.6){$\gamma$}
\put(3,.6){$\eta$}
\end{picture}
\end{center}

$B(d):=Q(3{\mathcal A})^{2,2}_1(d)$ is the quiver algebra of  $3{\mathcal A}$ modulo the relations
$$\beta\delta\eta=\beta\gamma\beta,\;\;\;\;
\delta\eta\gamma=\gamma\beta\gamma,\;\;\;\;
\eta\gamma\beta=d\cdot \eta\delta\eta,\;\;\;\;
\gamma\beta\delta=d\cdot \delta\eta\delta, \;\;\;\;\beta\delta\eta\delta=0,\;\;\;\;\eta\gamma\beta\gamma=0$$
for $d\in K\setminus\{0,1\}$.

Following \cite{Erdmann} the Cartan matrix of $B(d)$ is
$\left(\begin{array}{ccc}4&2&2\\ 2&3&1\\ 2&1&3\end{array}\right)$ and
the centre is $6$-dimensional. The Loewy series of the projective indecomposable modules
are given below.

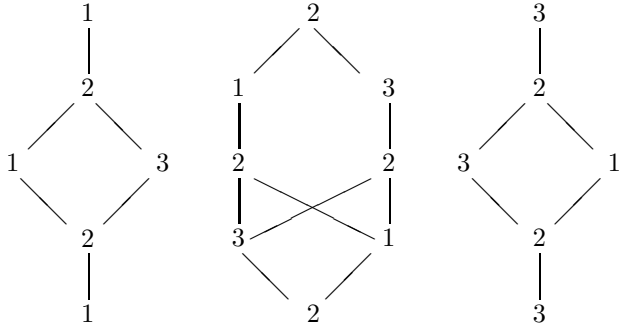
\unitlength1cm
\begin{center}
\begin{picture}(13,6)
\put(1,5){$1$}
\put(1,4){$2$}
\put(0,3){$1$}
\put(2,3){$3$}
\put(1,2){$2$}
\put(1,1){$1$}
\put(1.1,4.9){\line(0,-1){.6}}
\put(.9,3.9){\line(-1,-1){.6}}
\put(1.2,3.9){\line(1,-1){.6}}
\put(.2,2.9){\line(1,-1){.6}}
\put(1.9,2.9){\line(-1,-1){.6}}
\put(1.1,1.9){\line(0,-1){.6}}

\put(4,5){$2$}
\put(3,4){$1$}
\put(5,4){$3$}
\put(3,3){$2$}
\put(5,3){$2$}
\put(3,2){$3$}
\put(5,2){$1$}
\put(4,1){$2$}
\put(3.9,4.9){\line(-1,-1){.6}}
\put(4.1,4.9){\line(1,-1){.6}}
\put(3.1,3.9){\line(0,-1){.6}}
\put(5.1,3.9){\line(0,-1){.6}}
\put(3.1,2.9){\line(0,-1){.6}}
\put(5.1,2.9){\line(0,-1){.6}}
\put(3.3,2.9){\line(2,-1){1.6}}
\put(4.85,2.9){\line(-2,-1){1.6}}
\put(3.15,1.9){\line(1,-1){.6}}
\put(4.85,1.9){\line(-1,-1){.6}}

\put(7,5){$3$}
\put(7,4){$2$}
\put(6,3){$3$}
\put(8,3){$1$}
\put(7,2){$2$}
\put(7,1){$3$}
\put(7.1,4.9){\line(0,-1){.6}}
\put(6.9,3.9){\line(-1,-1){.6}}
\put(7.2,3.9){\line(1,-1){.6}}
\put(6.2,2.9){\line(1,-1){.6}}
\put(7.9,2.9){\line(-1,-1){.6}}
\put(7.1,1.9){\line(0,-1){.6}}

\end{picture}
\end{center}

We obtain a basis of the socle of $B(d)$ by
$$\{s_1:=\beta\delta\eta\gamma,s_2:=\eta\gamma\beta\delta,s_3:=\gamma\beta\delta\eta\}.$$
The closed paths of the algebra are
$$\{e_0,e_1,e_2,\beta\gamma, \gamma\beta, \delta\eta, \eta\delta,
\beta\delta\eta\gamma, \eta\gamma\beta\delta, \gamma\beta\delta\eta\}.$$
The centre is formed by linear combinations of closed paths and has a basis
$$\{1,\beta\gamma+\gamma\beta+\frac 1d\eta\delta,
\beta\gamma+\delta\eta+\eta\delta,
\beta\delta\eta\gamma,\eta\gamma\beta\delta,\gamma\beta\delta\eta\}$$
as is easily verified.
Non closed paths are clearly commutators. Obviously
$$\beta\delta\eta\gamma\equiv\eta\gamma\beta\delta\equiv\gamma\beta\delta\eta
\text{ mod }[B(d),B(d)].$$
Moreover,
$$\beta\gamma-\gamma\beta\in [B(d),B(d)]\mbox{ and }\delta\eta-\eta\delta\in [B(d),B(d)]. $$
Since the dimension of the centre of $B(d)$ coincides with the dimension of the commutator quotient, we get a basis of $B(d)/[B(d),B(d)]$
by $$\{e_0,e_1,e_2,\beta\gamma,\delta\eta,\beta\delta\eta\gamma\}.$$
We now suppose that the characteristic $p$ of $K$ is $p=2$. If $x$ is a square in $K$, then denote
$y=\sqrt x$ if $y^2=x$. We compute
$$
(\beta\gamma)^2=\beta\gamma\beta\gamma=\beta\delta\eta\gamma\mbox{ and }
(\delta\eta)^2=\delta\eta\delta\eta=\frac{1}{d}\delta\eta\gamma\beta.
$$
If $d$ is a square in $K$, then
$$((\gamma\beta)+\sqrt{d}(\eta\delta))^2=
\gamma\eta\gamma\beta+d\eta\delta\eta\delta=
\delta\eta\gamma\beta+\eta\gamma\beta\delta\in[B(d),B(d)]$$
so that $T_1(B(d))/[B(d),B(d)]$ is $1$-dimensional.
If $d$ is not a square, then $T_1(B(d))=[B(d),B(d)]$.

Let us consider the centre. Denote $\beta\gamma+\gamma\beta+\frac 1d\eta\delta=x$ and
$\beta\gamma+\delta\eta+\eta\delta=y$. Then we get
\begin{eqnarray*}
x^2=(\beta\gamma+\gamma\beta+\frac 1d\eta\delta)^2&=&
\beta\delta\eta\gamma+\frac 1d\delta\eta\gamma\beta+\frac 1d\eta\gamma\beta\delta
=s_1+\frac 1d s_2+\frac 1d s_2,\\
y^2=(\beta\gamma+\delta\eta+\eta\delta)^2&=&
\beta\delta\eta\gamma+\delta\eta\gamma\beta+\frac{1}{d^3}\eta\gamma\beta\delta
=s_1+s_2+\frac{1}{d^3}s_3,\\
xy=(\beta\gamma+\delta\eta+\eta\delta)\cdot(\beta\gamma+\gamma\beta+\frac 1d\eta\delta)&=&
\beta\delta\eta\gamma+\delta\eta\gamma\beta+\frac{1}{d^2}\eta\gamma\beta\delta
=s_1+s_2+\frac{1}{d^2}s_3.
\end{eqnarray*}
The coefficient matrix above has determinant $\frac{(d-1)^2}{d^4}$ and since $d\neq 1$,
the elements $x^2,y^2,xy$ are linearly independent, and hence
$Z(B(d))\simeq K[x,y]/(x^3,y^3,x^2y,xy^2)$. Moreover, choose the Frobenius
form given by
$$\psi(\beta\delta\eta\gamma)=\psi(\delta\eta\gamma\beta)=\psi(\eta\gamma\beta\delta)=1
\mbox{ and $\psi(c)=0$ if $c$ is a path of length at most $3$,}$$ following
Remark~\ref{symmetrisingform}. The relations are homogeneous, which shows that in order to prove symmetry of the form
we only need to consider paths $c_1$ and $c_2$ such that the lengths of $c_1$ and
$c_2$ sum up to $4$. The verification is a trivial and short computation which can be left to the reader.

Suppose now that $K$ is a perfect field. An elementary computation gives that $T_1^\perp(B(d))$
has a basis $\{x,s_1,s_2,s_3\}$, and therefore
$Z(B(d))/T_1^\perp(B(d))\simeq K[y]/y^2$, independently of $d$.

\begin{Theorem} \label{Qthreesimpletheorem}
Let $K$ be a field of characteristic $2$, and let $B(d)$ be the algebra
$Q(3{\mathcal A})^{2,2}_1(d)$. Then $\dim_K(T_1^\perp(B(d))/R(B(d)))=1$ if $d$
is a square in $K$, and $\dim_K(T_1^\perp(B(d))/R(B(d)))=0$ if $d$
is not a square in $K$.
\end{Theorem}

Proof: is done above.
\dickebox

\begin{Rem}
Unlike in case of Theorem~\ref{main} and its Corollary~\ref{Finalcorollarytwosimples}, 
using K\"ulshammer ideals we cannot distinguish the derived category of 
$Q(3{\mathcal A})^{2,2}_1(d)$ from the derived category of 
$Q(3{\mathcal A})^{2,2}_1(d')$ for two parameters $d,d'$.
If $K$ is perfect of characteristic $2$, then all elements of $K$ are squares. 
Theorem~\ref{prop:zimmermann} needs that $K$ is perfect for
the invariance of K\"ulshammer ideals under derived equivalences and
$K$ is even algebraically closed for the invariance under stable equivalences
of Morita type.
We can only say that the algebra $Q(3{\mathcal A})^{2,2}_1(d)$ is not isomorphic to the 
algebra $Q(3{\mathcal A})^{2,2}_1(d')$ if $d$ is a square and $d'$ is not. 
\end{Rem}

\end{document}